\documentclass[12pt,a4paper]{article}
\usepackage{amssymb}
\usepackage{amsmath}

\newtheorem{theorem}{Theorem}[section]
\newtheorem{definition}[theorem]{Definition}
\newtheorem{lemma}[theorem]{Lemma}

\date{}

\begin{document}

\title{Anti-commutative Gr\"{o}bner-Shirshov basis of a free Lie algebra}

\author{
L. A. Bokut\footnote {Supported by the RFBR  and
the Integration Grant of the SB RAS (No. 1.9).} \\
{\small \ School of Mathematical Sciences, South China Normal
University}\\
{\small Guangzhou 510631, P. R. China}\\
{\small Sobolev Institute of Mathematics, Russian Academy of
Sciences}\\
{\small Siberian Branch, Novosibirsk 630090, Russia}\\
{\small Email: bokut@math.nsc.ru}\\
\\
 Yuqun
Chen\footnote {Corresponding author.} \footnote{Supported by the
NNSF of China (No.10771077) and the NSF of Guangdong Province (No.06025062).} \  and Yu Li\\
{\small \ School of Mathematical Sciences, South China Normal
University}\\
{\small Guangzhou 510631, P. R. China}\\
{\small Email: yqchen@scnu.edu.cn}\\
{\small LiYu820615@126.com}}

\maketitle \noindent\textbf{Abstract:} One of the natural ways to
prove that the Hall words (Philip Hall, 1933) consist of a basis of
a free Lie algebra is a direct construction: to start with a linear
space spanned by Hall words, to define the Lie product of Hall
words, and then to check that the product yields the Lie identities
(Marshall Hall, 1950). Here we suggest another way using the
Composition-Diamond lemma for free anti-commutative
(non-associative) algebras (A.I. Shirshov, 1962).



\section{Introduction}

The history of the Hall basis (M. Hall, 1950,~\cite{MH}) of a free
Lie algebra is rather peculiar. P. Hall (1933,~\cite{PH}) invented
the basic commutators for a free group. From his paper, it followed
that any element of a free Lie algebra is a linear combination of
the basic Lie commutators. W. Magnus (1937,~\cite{M}) and E. Witt
(1937,~\cite{W}) proved that the Lie algebra obtained from a free
associative algebra is free. By the way, a simple proof of this
result had been lately found by A.I. Shirshov \cite{S2} (using
Lyndon--Shirshov words). Using this result, E. Witt found the Witt
formula for dimension of $L^n/L^{n-1}$, where L is a free Lie
algebra (the same paper). From the Witt formula, it may be seen that
basic Lie monomials must be linearly independent in the free Lie
algebra; according to K.W. Gruenberg \cite{G}, P. Hall knew this
fact. Despite all of that, M. Hall \cite{MH} was the first who
formulated and proved that ``the Hall words" constitute a linear
basis of a free Lie algebra. He used a direct construction of a free
Lie algebra: to start with a linear space spanned by Hall words, to
define the Lie product of Hall words, and then to check that the
product yields the Lie identities. For more general words, the
analogous method had been used lately by A.I. Shirshov \cite{S4}.
Shirshov's bases of a free Lie algebra, or better to say
Hall-Shirshov's bases (cf. \cite{R}), contain the Hall basis, the
Lyndon-Shirshov basis (\cite{S2}, \cite{L}), as well as bases that
lead to bases of free solvable (poly-nilpotent) Lie algebras
\cite{B}. By the way they do not contain the left normed basis by
E.S. Chibrikov \cite{Ch}.

In this paper, we are giving a detail proof of the Shirshov's
(Composition--Diamond) lemma for free anti-commutative algebras,
mostly following his original paper \cite{S3} with some improvement
in the terminology following Bokut \cite{b72}. (Shirhsov avoids
``composition of inclusion" proving that starting with  any finite
set of ``polynomials", one can effectively find an ``equivalent" set
with no composition of inclusion at all). As far as we understand,
this lemma was a step toward to a much more involved but of the same
kind of lemma for free Lie algebras (see \cite{Sh}). The last paper
is now wildly recognized as a pioneering paper in the theory of
G\"obner--Shirshov bases for Lie and associative algebras. Let us
recall that the same kind of theory for commutative associative
algebras is mostly due to B. Buchberger \cite{bu70} though some
ideas were discovered by H. Hironaka \cite{HH}. Speaking about
Shirshov's lemma for free anti-commutative algebras, it is of course
a (main) part of his ``Gr\"obner--Shirshov bases" theory for
anti-commutative algebras from the same paper \cite{S3}. His main
application was a simple and conceptual proof of the algorithmic
decidability of the word problem for any finitely presented
anti-commutative (commutative) algebra. It is an analogy of Zhukov's
theorem for non-associative algebras (cf. \cite{Z}).

 We present here also a new proof that the Hall words are linearly independent
 based on Shirshov's Composition-Diamond lemma for anti-commutative algebras
\cite{S3}. We use also a direct construction of a free
anti-commutative algebra following A.I. Shirshov \cite{S1}.

\section{Direct construction of a free anti-commutative algebra AC(X)}

Let $X=\{x_i|i\in I \}$ be a well ordered set, $X^*$ the set of all
associative words $u$ in $X$, and $X^{\ast \ast }$ the set of all
non-associative words $(u)$ in $X$. We assume $(u)$ as a bracketing
of $u$. As a rule, we will omit ``non-associative" in the expression
``non-associative word".
     Then we define normal words $N=\{[u]\}$ and order of them by induction
 on the length $n=|[u]|$ of $[u]$:\\

$(i)$ If $n=1$, then $[u]=x_{i}$ is a normal word. Define $%
x_{i}>x_{j}$ if $i>j$.\\

Let $N_{n-1}=\{[u]|[u]\ is\ a\ normal\ word\ and\ |[u]|\leq
{n-1}\}$, $n>1$ and suppose that $``<"$ is a total order on
$N_{n-1}$. Then

\begin{enumerate}
\item[(ii)] \ If $n>1$ and $(u)=((v)(w))$ is a word of length $n$, then $(u)$
is a normal word, if and only if
\begin{enumerate}
\item[(a)] \ both $(v)$ and $(w)$ are normal words, that is, $(v)=[v]$ and $(w)=[w]$,
and

\item[(b)] \ $[v]>[w]$.
\end{enumerate}
\end{enumerate}
 Define $``<"$: Let $[u]$, $[v]$ be normal
words of length $\leq n$. Then $[u]<[v]$, if and only if one of the
following three cases holds:
\begin{enumerate}
\item[(a)] \ $|[u]|<n$, $|[v]|<n$
and $[u]<[v]$.

\item[(b)] \ $|[u]|<n$ and $|[v]|=n$.

\item[(c)] \ If $|[u]|=|[v]|=n$, $[u]=[[u_{1}][u_{2}]]$ and $[v]=[[v_{1}][v_{2}]]$, then $%
[u_{1}]<[v_{1}]$ or $([u_{1}]=[v_{1}]\ and\ [u_{2}]<[v_{2}])$.
\end{enumerate}

It is clear that the order ``$<$" on $N$ is a well order. This order
is called deg-lex (degree-lexicographical) order and we use this
order through this paper.

Let $k$ be a field and $AC(X)$ be a $k$-space spanned by normal
words. Let us define the product of normal words by the following
way:\newline
\begin{equation*}
\lbrack u][v]=\left\{
\begin{array}{r@{\quad:\quad}l}
[[u][v]] & [u]>[v] \\
-[[v][u]] & [u]<[v] \\
0\text{ \ \ } & [u]=[v]%
\end{array}%
\right.
\end{equation*}

\noindent \textbf{Remark} \ By definition, for any $(u)\in X^{**}$,
there exists a unique $[v]\in N$ such that $(u)=\pm [v]$ or 0. We
will denote $[v]$ by $\widetilde{(u)}$ sometimes if $(u)\neq0$.

\begin{theorem}
(A.I. Shirshov,~\cite{S1})~$AC(X)$ is a free anti-commutative
algebra generated by $X$.
\end{theorem}
{\bf Proof.} Suppose that $ f=\sum\limits_i \alpha_i[u_i], \
g=\sum\limits_j\beta_j[v_j]\in AC(X),$ where $\alpha_i,\beta_j\in k$
and $[u_i], [v_j]$ are normal  words. Then
$$
fg+gf=\sum\limits_{i,j}\alpha_i \beta_j([u_i][v_j]+[v_j][u_i])=0.
$$
So, $AC(X)$ is an anti-commutative algebra. We will prove $AC(X)$ is
free on $X$. Let $A$ be an anti-commutative algebra and $\xi$ be a
map from $X$ to $A$. Then we define
$$
\bar \xi: AC(X)\longrightarrow A, \ \  [x_{i_1}x_{i_2}\cdots
x_{i_n}]\mapsto (\xi{(x_{i_1})}\xi{(x_{i_2})}\cdots \xi{(x_{i_n})}),
$$
 where $[x_{i_1}x_{i_2}\cdots
x_{i_n}]$ is a normal word. It is easy to check that $\bar \xi$ is
the unique algebra homomorphism such that
$$
\bar{\xi}\circ i=\xi,
$$
where $i: \ X\rightarrow AC(X)$ is the including mapping. This
completes our proof. \ \ \ \ $\square$

\section{Composition-Diamond Lemma for $AC(X)$}

In this section, we establish Composition-Diamond lemma for the free
anti-commutative algebra $AC(X)$.

Let $N=\{[u]|~[u]$\ is\ a\ normal\ word $\}$ and $``<"$ be the
deg-lex order on $N$ as before. It is easy to check that $``<"$ is a
monomial well order in the following sense:
\begin{equation}\label{e1}
[u]>[v]\Rightarrow \widetilde{([u][w])}>\widetilde{([v][w])}
\end{equation}
where $[w]\ne [u],\ [w]\ne [v]\ $. As a result, we have:
\begin{equation}\label{e2}
[u]>[v]\Rightarrow [a[u]b]>\widetilde{(a[v]b)}
\end{equation}
where $[a[u]b]$ is a normal word with subword $[u]$ and
$(a[v]b)=[a[u]b]|_{[u]\mapsto [v] }$.

Given a polynomial $f\in AC(X)$, it has the leading word
$[\overline{f}] \in N$ according to the above order on $N$, such
that
$$
f=\sum_{[u]\in N}f([u])[u]=\alpha
[\overline{f}]+\sum{\alpha}_i[u_i],
$$
where $[\overline{f}]>[u_i], \ \alpha , \ {\alpha}_i, \ f([u])\in
k$. We call $[\overline{f}]$ the leading term of $f$. Denote the set
$\{ [u]|f([u])\neq 0\}$  by $suppf$ and $deg(f)$ by
$|\overline{f}|$. $f$ is called monic if $\alpha=1$.

\begin{definition}
Let $S\subset AC(X)$ be a set of monic polynomials, $s\in S$ and $(u)\in X^{**}$. We define $S$-word $%
(u)_s$ by induction:
\begin{enumerate}
\item[(i)] $(s)_s=s$ is an $S$-word of $S$-length
1.
\item[(ii)] If $(u)_s$ is an $S$-word of $S$-length k and $(v)$
is a nonassociative word of length $l$, then
\begin{equation*}
(u)_s(v)\ and\ (v)(u)_s
\end{equation*}
are $S$-words of length $k+l$.
\end{enumerate} The $S$-length of an
$S$-word $(u)_s$ will be denoted by $|u|_s$.
\end{definition}

\begin{definition}
$S$-word $(u)_s$ is called normal $S$-word, if $(u)_{[\bar s]
}=(a[\bar s ]b)$ is a normal word. We denote $(u)_s$ by $[u]_s$, if
$(u)_s$ is a normal $S$-word. We also call the normal $S$-word
$[u]_s$ to be normal $s$-word. From (\ref{e2}) it follows that
$\overline{[u]_s}=[u]_{[\bar s] }$.
\end{definition}

Let $f,g$ be monic polynomials in $AC(X)$. Suppose that there exist
$a,b\in {X^*} $ such that $[\bar f] =[a[\bar g] b]$, where $[agb]$
is a normal $g$-word. Then we set $[w]=[\bar f ]$ and define the
composition of inclusion
\begin{equation*}
(f,g)_{[w]}=f-[agb].
\end{equation*}
We note that
\begin{equation*}
(f,g)_{[w]}\in Id(f,g)\ \ and\ \ \overline{(f,g)_{[w]}}<[w].
\end{equation*}
Transformation, $f\longmapsto f-[agb]$ is called the Elimination of
Leading Word (ELW) of $g$ in $f$.

Given a nonempty subset $S\subset AC(X)$, we shall say that the
composition $(f,g)_{[w]}$ is trivial modulo $(S,[w])$, if
\begin{equation*}
(f,g)_{[w]}=\sum\limits_i\alpha_i[a_is_ib_i],
\end{equation*}
where each $\alpha_i\in k,\ a_i,b_i\in X^*,\ s_i\in S,\
[a_is_ib_i]$ is normal S-word and $[a_i[\bar{s_i}]b_i]<{[w]}$. If
this is the case, then we write $(f,g)_{[w]}\equiv 0\ mod
(S,{[w]})$.

Let us note that if $(f,g)_{[w]}$ goes to $0$ by ELW's of $S$,
then $(f,g)_{[w]}\equiv 0$ mod$(S,[w])$. Indeed, using ELW's of
$S$, we have
$$
(f,g)_{[w]}\longmapsto (f,g)_{[w]}-\alpha_1[a_1s_1b_1]=f_2
\longmapsto f_2-\alpha_2[a_2s_2b_2]\longmapsto \cdots \longmapsto
0.
$$
So, $(f,g)_{[w]}=\sum\limits_i\alpha_i[a_is_ib_i]$ where
$[a_i[\bar{s_i}]b_i] \leq \overline{(f,g)_{[w]}}<{[w]}$.

In general, for $p,q\in AC(X)$, we write
$$
p\equiv q\quad mod(S,[w])
$$
which means that $p-q=\sum\alpha_i [a_i s_i b_i] $, where each
$\alpha_i\in k, \ a_i,b_i\in X^{*}, \ s_i\in S$ and $[a_i [\bar
{s_i}] b_i]<[w]$.

\begin{definition}
Let $S\subset AC(X)$ be a nonempty set of monic polynomials and the
order ``$<$" as before. Then the set $S$ is called a
Gr\"{o}bner-Shirshov basis, if any composition $(f,g)_{[w]}$ with
$f,g\in S$ is trivial modulo $(S,{[w]})$, i.e., $(f,g)_{[w]}\equiv0$
$mod(S,{[w]})$.
\end{definition}

\begin{lemma}\label{3.4}
Let $[v]_s$ be a normal $s$-word and $[w]$, $[w']$ normal words. If
$[w]>[w']$ and $[w]>[v]_{[\bar s]}$, then $s$-word $([w'][v]_s)$ has
a representation:
$$
([w'][v]_s)=-([v]_s[w'])=\sum\limits_i\alpha_i[u_i]_s,
$$
where each $\alpha_i\in k$, $[u_i]_s$ normal $s$-word and
$[u_i]_{[\bar s]}<min\{[[w][v]_{[\bar s]}],[[w][w']]\}$.
\end{lemma}

\noindent\textbf{Proof.}  Suppose that
$$
[v]_s=\gamma[v]_{[\bar s]}+\sum\limits_n \gamma_n[v_n],
$$
where $0\ne \gamma,\gamma_n\in k,\ [v]_{[\bar s]},[v_n]\in N\ and\
[v_n]<[v]_{[\bar s]}$. Now we consider the following three cases:
\begin{enumerate}
\item[(a)] If $[v]_{[\bar s]}<[w']$, then
$([w'][v]_s)$ is already a normal $s$-word and $[[w'][v]_{[\bar
s]}]<min\{[[w][v]_{[\bar s]}],[[w][w']]\}$.
\item[(b)] If $%
[v]_{[\bar s]}>[w']$, then
\begin{equation*}
([w'][v]_s)=-([v]_s[w']).
\end{equation*}
Here $([v]_s[w'])$ is a normal $s$-word and $[[v]_{[\bar s]}[w']]<min\{[[w][v]_{[\bar s]}],[[w][w']]\}$.\\
\item[(c)] If $[v]_{[\bar s]}=[w']$, then
\begin{equation*}
([w'][v]_s)=([v]_{[\bar s]}[v]_s)=\sum\limits_n
\gamma^{-1}\gamma_n([v]_s[v_n]),
\end{equation*}
since
\begin{equation*}
0=[v]_s[v]_s=(\gamma[v]_{[\bar s]}+\sum\limits_n
\gamma_n[v_n])[v]_s.
\end{equation*}
Now, clearly, each $([v]_s[v_n])$ is normal $s$-word and
$[[v]_{[\bar s]}[v_n]]<min\{[[w][v]_{[\bar s]}],[[w][w']]\}$. \ \ \
\ $\square$
\end{enumerate}

\begin{lemma}\label{3.5}
Let $(u)_s$ be an $S$-word. Then  $(u)_s$ has a representation:
$$
(u)_s=\sum\limits_i \alpha_i [u_i]_{s},
$$
where each $\alpha_i\in k$ and $[u_i]_{s}$ is normal $s$-word.
\end{lemma}

\noindent\textbf{Proof.} We use induction on $|u|_s$. If $|u|_s=1$, then $%
(u)_s=s$ and the result holds. If $|u|_s>1$, then $(u)_s=(v)_s(w)$ or $%
(u)_s=(w)(v)_s$. Here we consider the case $(u)_s=(v)_s(w)$. The
other one is similarly proved. By induction,
\begin{equation*}
(v)_s=\sum\limits_j\beta_j[v_j]_{s},
\end{equation*}
where $\beta_j\in k \ and \ [v_j]_{s}$ is normal $s$-word. Without
loss of generality, we may assume $(v)_s$ is a normal $s$-word and
$(w)$ is a normal word. Then $(u)_s=([v]_s[w])$. Just like the proof
in the Lemma \ref{3.4}, we know that $(u)_s=([v]_s[w])$ is a linear
combination of normal $s$-words.

The proof is completed.   \ \ $\square$

\begin{lemma}\label{a1}
Let $[u]_s=[a'[v]a''sb]$ (or $[u]_s=[asb'[v]b'']$) be a normal
$s$-word, $[w]$ a normal word and $[v]>[w]$. Then, the $s$-word
$(u')_s=(a'[w]a''sb)=[a'[v]a''sb]|_{[v]\mapsto [w]}$ has a
representation:
\begin{equation} \label{a}
(u')_s =\sum_i\alpha_i[u_i]_s,
\end{equation}
where each $\alpha_i\in k$, $[u_i]_s$ normal $s$-words and
$[u_i]_{[\bar s]}<[u]_{[\bar s]}$.

For the $s$-word $(u')_s=(asb'[w]b'')=[asb'[v]b'']|_{[v]\mapsto
[w]}$, it has a similar representation to (\ref{a}).
\end{lemma}
\noindent\textbf{Proof.} We prove only the first case. The other one
is similarly proved. Induction on $|u|_s$. If $|u|_s=|v|+1$, then
$[u]_s=[[v]s]$ and $(u')_s=([w]s)$.  Then the result follows from
Lemma \ref{3.4}. Suppose that $|u|_s>|v|+1$ and
$[u]_s=[[u_1][u_2]_s]$ or $[u]_s=[[u_1]_s[u_2]]$. We deal with only
 the case $[u]_s=[[u_1][u_2]_s]$. If $[v]$ is a subword of
$[u_1]$, then we let $(u_1^*)=[u_1]|_{[v]\mapsto[w]}$ and let
$(u_1^*)=\widetilde{(u_1^*)}$ (if $(u_1^*)=0$ the case is trivial).
Since $[v]>[w]$, by (\ref{e2}), we have $[u_1]>\widetilde{(u_1^*)}$.
Now, $(u')_s=(\widetilde{(u_1^*)}[u_2]_s)$ and the result follows
from Lemma \ref{3.4}. If $[v]$ is a subword of $[u_2]_s$, then by
induction we have $(u_2')_s=\sum\alpha_i[u_i]_s$, where each
$[u_i]_{[\bar s]}<[u_2]_{[\bar s]}$. Then
$(u')_s=([u_1](u_2')_s)=\sum\alpha_i[u_1][u_i]_s$. So, by Lemma
\ref{3.4} again, we get the result.        \ \ $\square$

\begin{lemma}\label{a2}
Let $[u]_s=[asb]$, $[v]_t$ be normal $s$- and $t$- words
respectively. If $[v]_{[\bar t]}<[\bar s]$, then the $t$-word
$(u)_{[v]_t}=[asb]|_{s\mapsto [v]_t}$ has a representation:
$$
(u)_{[v]_t}=\sum_i\alpha_i[u_i]_t,
$$
where each $\alpha_i\in k$, $[u_i]_t$ normal $t$-words and
$[u_i]_{[\bar t]}<[u]_{[\bar s]}$.
\end{lemma}

\noindent\textbf{Proof.} By induction on $|u|_s$ and Lemma
\ref{3.4}, we may easily get the result. \ \ $\square$

\begin{lemma}\label{3.7}
Let $[a_1s_1b_1],~[a_2s_2b_2]$ be normal $S$-words. If $S$ is a
Gr\"{o}bner-Shirshov basis in $AC(X)$ and
$[w]=[a_1[\overline{s_1}]b_1]=[a_2[\overline{s_2}]b_2]$, then
\begin{equation*}
[a_1s_1b_1]\equiv [a_2s_2b_2] \ \ mod (S,[w]).
\end{equation*}

\end{lemma}

\noindent\textbf{Proof.} We have $a_1\bar s_1 b_1=a_2\bar s_2 b_2$
as associative words in the alphabet $X\bigcup \{\bar {s_1}, \bar
{s_2}\}$. There are two cases to consider.

Case 1. Suppose that subwords $\bar s_1$ and $\bar s_2$ of $w$ are
disjoint, say, $|a_2|\geq |a_1|+|\bar s_1|$. Then, we can assume
that
$$
a_2=a_1\bar s_1 c \ \ and  \ b_1=c\bar s_2 b_2
$$
for some $c\in X^*$, and so, $ [w]=[a_1[\bar s_1] c [\bar s_2] b_2].
$ Now,
\begin{eqnarray*}
[a_1 s_1 b_1]-[a_2 s_2 b_2]&=&[a_1 s_1 c [\bar s_2]
b_2]-[a_1[\bar s_1 ]c  s_2 b_2]\\
&=&[a_1 s_1 c [\bar s_2] b_2]-(a_1 s_1 c  s_2 b_2)+(a_1 s_1 c  s_2
b_2)-[a_1[\bar s_1] c  s_2 b_2]\\
&=&(a_1 s_1 c ([\bar s_2] - s_2) b_2)+(a_1(s_1-[\bar s_1]) c  s_2
b_2).
\end{eqnarray*}

Since $[\overline{[\overline{s_2}]-s_2}]<[\bar s_2]$ and
$[\overline{s_1-[\overline{s_1}]}]<[\bar s_1]$, and by the Lemma
\ref{a1}, we conclude that
$$
[a_1 s_1 b_1]-[a_2 s_2 b_2]=\sum\limits_i
\alpha_i[u_is_1v_i]+\sum\limits_j \beta_j[u_js_2v_j]
$$
for some $\alpha_i,\beta_j\in k$, normal S-words $[u_is_1v_i]$ and
$[u_js_2v_j]$ such that $ [u_i[\bar s_1]v_i],[u_j[\bar s_2]v_j]<[w].
$ So,
$$
[a_1s_1b_1]\equiv [a_2s_2b_2]\ mod (S,[w]).
$$

Case 2. Suppose that the subword $\bar s_1$ of $w$ contains $\bar
s_2$ as a subword. We assume that
$$
[\bar s_1]=[a[\bar s_2]b], \ a_2=a_1a  \mbox{ and } b_2=bb_1, \mbox{
that is, } [w]=[a_1[a[\bar s_2]b]b_1]
$$
for the normal $S$-word $[a s_2 b]$. We have
\begin{eqnarray*}
[a_1 s_1 b_1]-[a_2 s_2 b_2]&=&[a_1 s_1 b_1]-[a_1 [a s_2 b] b_1]\\
&=&(a_1(s_1-[as_2b])b_1)\\
&=&(a_1(s_1,s_2)_{[w_1]}b_1),
\end{eqnarray*}
where $[w_1]=[\overline{s_1}]=[a[\bar s_2] b]$. Since $S$ is a
Gr\"{o}bner-Shirshov basis,
$(s_1,s_2)_{[w_1]}=\sum\limits_i\alpha_i[c_i s_i d_i]$ for some
$\alpha_i\in k$, normal S-words $[c_is_id_i]$ with each $[c_i[\bar
s_i] d_i]<[w_1]=[\bar{s_1}]$. By Lemma \ref{a2}, we have
\begin{eqnarray*}
&&[a_1 s_1 b_1]-[a_2 s_2 b_2]=(a_1(s_1,s_2)_{[w_1]}b_1)\\
&=&\sum\limits_i\alpha_i(a_1[c_is_id_i]b_1)=\sum\limits_j\beta_j[a_js_jb_j]
\end{eqnarray*}
for some $\beta_j\in k $, normal S-words $[a_js_jb_j]$ with each
$[a_j[\bar s_j ]b_j]<[w]=[a_1[\bar {s_1}]b_1]$. \\
So,
$$
[a_1s_1b_1]\equiv [a_2s_2b_2]\ mod (S,[w]).  \ \ \ \ \square
$$

\begin{lemma}\label{3.8}
Let $S\subset AC(X)$ be set of monic polynomials and
$Red(S)=\{[u]\in N |[u]\ne [a[\bar s] b]\ a,b\in X^*,\ s\in S \mbox{
and } [as b] \mbox{ is a normal } S\mbox{-word}\}$. Then for any
$f\in AC(X)$,
\begin{equation*}
f=\sum\limits_{[u_i]\leq [\bar f] }\alpha_i[u_i]+
\sum\limits_{[a_j[\overline{s_j}]b_j]\leq[\bar
f]}\beta_j[a_js_jb_j],
\end{equation*}
where each $\alpha_i,\beta_j\in k, \ [u_i]\in Red(S)$ and
$[a_js_jb_j]$ normal $S$-word.
\end{lemma}
{\bf Proof.} Let  $f=\sum\limits_{i}\alpha_{i}[u_{i}]\in{AC(X)}$,
where $0\neq{\alpha_{i}\in{k}}$ and $[u_{1}]>[u_{2}]>\cdots$. If
$[u_1]\in{Red(S)}$, then let $f_{1}=f-\alpha_{1}[u_1]$. If
$[u_1]\not\in{Red(S)}$, then there exist some $s\in{S}$ and
$a_1,b_1\in{X^*}$, such that $[\bar f]=[a_1[\bar{s_1}]b_1]$. Let
$f_1=f-\alpha_1[a_1s_1b_1]$. In both cases, we have
$[\bar{f_1}]<[\bar{f}]$. Then the result follows from the induction
on $[\bar{f}]$. \ \ \ \ $\square$

\begin{theorem}(Shirshov \cite{S3})
Let $S\subset AC(X)$ be a nonempty set of monic polynomials and the
order $``<"$ as before. Then the following statements are
equivalent:
\begin{enumerate}
\item [(i)] $S$ is a Gr\"{o}bner-Shirshov basis.

\item [(ii)] $f\in Id(S)\Rightarrow [\bar f] =[a[\bar s ]b]$ for some $s\in S\
and\ a,b\in X^*$, where $[as b]$ is  normal $S$-word.

\item [(ii)'] $f\in Id(S)\Rightarrow  f = \alpha_1[a_1s_1
b_1]+\alpha_2[a_2s_2b_2]+\cdots$, where $\alpha_i\in k, \
[a_1[\overline{s_1}]b_1]>[a_2[\overline{s_2}]b_2]>\cdots$ and each
$[as_i b]$ is  normal $S$-word.

\item [(iii)] $Red(S)=\{[u]\in N |[u]\ne [a[\bar s] b]\ a,b\in X^*,\ s\in S \mbox{ and }
[as b] \mbox{ is a normal } S\mbox{-word}\}$ is a basis of the
algebra $AC(X|S)$.
\end{enumerate}
\end{theorem}
{\bf Proof.} $(i)\Rightarrow (ii)$. \ Let $S$ be a
Gr\"{o}bner-Shirshov basis and $0\neq f\in Id(S)$. We can assume, by
Lemma \ref{3.5}, that
$$
f=\sum_{i=1}^n\alpha_i[a_is_ib_i],
$$
where each $\alpha_i\in k, \ a_i,b_i\in {X^*}, \ s_i\in S$ and $
[a_is_ib_i]$ normal $S$-word. Let
$$
[w_i]=[a_i[\overline{s_i}]b_i], \
[w_1]=[w_2]=\cdots=[w_l]>[w_{l+1}]\geq\cdots
$$
We will use the induction on $l$ and $[w_1]$ to prove that
$[\overline{f}]=[a[\overline{s}]b]$ for some $s\in S \ \mbox{and} \
a,b\in {X^*}$.

If $l=1$, then
$[\overline{f}]=\overline{[a_1s_1b_1]}=[a_1[\overline{s_1}]b_1]$ and
hence the result holds. Assume that $l\geq 2$. Then, by Lemma
\ref{3.7}, we have
$$
[a_1s_1b_1]\equiv[a_2s_2b_2] \ \ mod(S,[w]).
$$
Thus, if $\alpha_1+\alpha_2\neq 0$ or $l>2$, then the result holds.
For the case $\alpha_1+\alpha_2= 0$ and $l=2$, we use the induction
on $[w_1]$. Now, the result follows.\\

$(ii)\Rightarrow (ii)'$. \ Assume (ii) and $0\neq f\in Id(S)$. Let
$f=\alpha_1[\overline{f}]+\cdots$. Then, by (ii),
$[\overline{f}]=[a_1[\overline{s_1}]b_1]$. Therefore,
$$
f_1=f-\alpha_1[a_1s_1b_1], \ [\overline{f_1}]<[\overline{f}], \
f_1\in Id(S).
$$
Now, by using induction on $[\overline{f}]$, we have $(ii)'$.\\

$(ii)'\Rightarrow (ii)$. This part is clear.\\

$(ii)\Rightarrow(iii)$. Suppose that
$\sum\limits_{i}\alpha_i[u_i]=0$ in $AC(X|S)$, where $\alpha_i\in
k$, $[u_i]\in {Red(S)}$. It means that
$\sum\limits_{i}\alpha_i[u_i]\in{Id(S)}$ in ${AC(X)}$. Then all
$\alpha_i$ must be equal to zero. Otherwise,
$\overline{\sum\limits_{i}\alpha_i[u_i]}=[u_j]\in{Red(S)}$ for some
$j$ which contradicts (ii).

Now, for any $f\in{AC(X)}$, by Lemma \ref{3.8}, we have
$$
f=\sum\limits_{[u_i]\in Red(S),\ [u_i]\leq [\bar f] }\alpha_i[u_i]+
\sum\limits_{[a_js_jb_j]-normal,\ [a_j[\overline{s_j}]b_j]\leq[\bar
f]}\beta_j[a_js_jb_j].
$$
So, (iii) follows.\\

$(iii)\Rightarrow(i)$. For any $f,g\in{S}$ , by Lemma \ref{3.8}, we
have
$$
(f,g)_{[w]}=\sum\limits_{[u_i]\in Red(S),\ [u_i]<[w] }\alpha_i[u_i]+
\sum\limits_{[a_js_jb_j]-normal,\
[a_j[\overline{s_j}]b_j]<[w]}\beta_j[a_js_jb_j].
$$
Since $(f,g)_{[w]}\in {Id(S)}$ and by (iii), we have
$$
(f,g)_{[w]}= \sum\limits_{[a_js_jb_j]-normal,\
[a_j[\overline{s_j}]b_j]<[w]}\beta_j[a_js_jb_j].
$$
Therefore, $S$ is a Gr\"{o}bner-Shirshov basis. \ \ \ \ $\square$

\section{Gr\"{o}bner-Shirshov basis for a free Lie algebra}

In this section, we represent the free Lie algebra by the free
anti-commutative algebra and give a Gr\"{o}bner-Shirshov basis for a
free Lie algebra.

The proof of the following theorem is straightforward and we omit
the detail.

\begin{theorem}
Let $AC(X)$ be free anti-commutative algebra and let
\begin{equation*}
S=\{([u][v])[w]-([u][w])[v]-[u]([v][w]) \ | \ [u],[v],[w]\in N \
and\ [u]>[v]>[w]\}.
\end{equation*}
Then the algebra $AC(X|S)$ is the free Lie algebra generated by $X$.
\ \ \ \ $\square$
\end{theorem}

We now cite the definition of Hall words by induction on length:
\begin{enumerate}
\item[1)]$x_i$ is a Hall word for any $x_i\in X$.

Suppose we define Hall words of length $<n$.

\item[2)] Normal word $[[v][w]]$ is called Hall word if and
only if \begin{enumerate}

\item[(a)] both $[v]$ and $[w]$ are Hall words,

\item[(b)] if $[v]=[[v_1][v_2]]$, then $[v_2]\leq[w]$.
\end{enumerate}
\end{enumerate}

We denote $[u]$ by $[[u]]$, if $[u]$ is a Hall word. Let
\begin{eqnarray*}
 S_0&=&\{([[u]][[v]])[[w]]-([[u]][[w]])[[v]]-[[u]]([[v]][[w]]) \ |
 \\
&&[[u]]>[[v]]>[[w]]\  and  \ [[u]],[[v]],[[w]]\ are\ Hall\ words\}.
\end{eqnarray*}

\begin{lemma}\label{4.1}
Let $H$ be the set consisting of all Hall words. Then
\begin{equation*}
Red(S_0)=\{[u]\in N |[u]\ne [a[\bar s] b]\ a,b\in X^*,\ s\in S_0
\mbox{ and } [as b] \mbox{ is a normal } s\mbox{-word}\}=H.
\end{equation*}
\end{lemma}

\noindent\textbf{Proof.} Suppose $[u]\in Red(S_0)$. We will show
that $[u]$ is a Hall word by induction on $|[u]|=n$. If $n=1$, then
$$
[u]=x_i
$$
which is already a Hall word. Let $n>1$ and $[u]=[[v][w]]$. This
case has two subcases. By induction, we have that $ [v],[w] $ are
Hall words.\\

Subcase 1. If $|v|=1$, then $[u]$ is a Hall word.\\

Subcase 2. If $|[v]|>1\ and \ [v]=[[v_1][v_2]]$, then
\begin{equation*}
[v_2]\leq [w]
\end{equation*}
for $[u]\in Red(S_0)$. So, $[u]$ is a Hall word.

It's clear that every Hall word is in $Red(S_0)$ since every subword
of Hall word is also a Hall word.\ \ \ \ $\square$

The following lemma follows from Lemma \ref{3.8} and Lemma
\ref{4.1}.

\begin{lemma}\label{4.2}
In $AC(X)$, any normal word $[u]$ has the following presentation:
$$
[u]=\sum\limits_{i}\alpha_i[[u_i]]+\sum\limits_{j}\alpha'_j[u'_j]_{s'_j}
$$
 where $\alpha_i, \alpha'_j\in k$, $[[u_i]]$ are Hall words,
$[u'_j]_{s'_j}$ normal $S_0$-words, $s'_j\in S_0, \ [[u_i]],
[u'_j]_{[\overline {s'_j}]}\leq[u]$. Moreover, each $[[u_i]]$ has
the same length as $[u]$.
\end{lemma}

\begin{lemma}
Suppose $S$ and $S_0$ are sets defined as before. Then, in $AC(X)$,
we have
\begin{equation*}
Id(S)=Id(S_0).
\end{equation*}

\end{lemma}
{\bf Proof.} Since $S_0$ is a subset of $S$, it suffices to prove
that $AC(X|S_0)$ is a Lie algebra. We need only to prove that, in
$AC(X|S_0)$,
\begin{equation*}
([u][v])[w]-([u][w])[v]-[u]([v][w])=0,
\end{equation*}
where $[u],[v],[w]\in N$ and $[u]>[v]>[w]$. By Lemma \ref{4.2}, it
suffices to prove
\begin{equation*}
([[u]][[v]])[[w]]-([[u]][[w]])[[v]]-[[u]]([[v]][[w]])=0,
\end{equation*}
where  $[[u]]>[[v]]>[[w]]$. This is trivial by the definition of
$S_0$. \ \ \ \ $\square$

\begin{theorem}
Let the order $``<"$ be as before and
\begin{center}
 $S_0=\{([[u]][[v]])[[w]]-([[u]][[w]])[[v]]-[[u]]([[v]][[w]]) \ |
~[[u]]>[[v]]>[[w]]\  and  \ [[u]],[[v]],[[w]]\ are\ Hall\ words\}.$
\end{center}
 Then $S_0$ is a Gr\"obner-Shirshov basis in
$AC(X)$.
\end{theorem}
{\bf Proof.} To simplify notations, we write $u$ for $[[u]]$ and
$u_1u_2\cdots u_n$ for
$((((u_1)u_2)\cdots)u_n)$.\\
Let
\begin{equation*}
f_{uvw}=uvw-uwv-u(vw),
\end{equation*}
where $u,v,w $ are Hall words and $u>v>w$. It is easy to check that
$\overline{f_{uvw}}=uvw$.

Suppose $\overline{f_{u_1v_1w_1}}$ is a subword of
$\overline{f_{uvw}}$. Since $u,v,w$ are Hall words, we have
$u_1v_1w_1=uv,u=u_1v_1\ and \ v=w_1$. We will prove that the
composition
\begin{equation*}
(f_{uvw},f_{u_1v_1w_1})_{uvw}
\end{equation*}
is trivial modulo $(S_0, uvw)$. We note that $u_1>v_1>w_1=v>w$.

Firstly, we  prove that the following statements hold
mod$(S_0,uvw)$:
\begin{enumerate}

\item[1)]$u_1vv_1w-u_1vwv_1-u_1v(v_1w)\equiv 0$.

\item[2)]$u_1(v_1v)w-u_1w(v_1v)-u_1(v_1vw)\equiv 0$.

\item[3)]$u_1wv_1v-u_1wvv_1-u_1w(v_1v)\equiv 0$.

\item[4)]$u_1(v_1w)v-u_1v(v_1w)-u_1(v_1wv)\equiv 0$.

\item[5)]$u_1(vw)v_1-u_1(vwv_1)-u_1v_1(vw)\equiv 0$.

\item[6)]$u_1vwv_1-u_1wvv_1-u_1(vw)v_1\equiv 0$.

\item[7)]$u_1v_1wv-u_1wv_1v-u_1(v_1w)v\equiv 0 $.

\item[8)]$u_1(v_1vw)-u_1(v_1wv)-u_1(v_1(vw))\equiv 0$.

\end{enumerate}

We only prove $5)$. $1)$--$4)$ are similarly proved to $5)$ and
$6)$--$8)$ follow from ELW's of $S_0$.

By ELW's of $S_0$, we may assume, without loss of generality, that
$vw$ is a Hall word. It's easy to check 5) holds in the following
three cases: $vw>v_1,\ vw=v_1$ and $vw<v_1$. For example, let
$vw>v_1$ and $E=u_1(vw)v_1-u_1(vwv_1)-u_1v_1(vw)$. We consider the
following cases: if $u_1>vw$, then $E=f_{u_1(vw)v_1}\equiv 0$; if
$u_1=vw$, then $E=0$; if $u_1<vw$, then $E=-f_{(vw)u_1v_1}\equiv 0$.
So $5)$ is proved.

Secondly, we have
\begin{align*}
(f_{uvw},f_{u_1v_1w_1})_{uvw}&=f_{uvw}-(f_{u_1v_1w_1})w \\
&=u_1vv_1w+u_1(v_1v)w-u_1v_1wv-u_1v_1(vw).
\end{align*}
Let
\begin{equation*}
A=u_1vv_1w+u_1(v_1v)w \ \ and \ \ B=-u_1v_1wv-u_1v_1(vw).
\end{equation*}

Then, by $1)$--$8)$, we have
\begin{align*}
A & \equiv u_1vwv_1+u_1v(v_1w)+u_1w(v_1v)+u_1(v_1vw) \ \ \ \\
&\equiv
u_1wvv_1+u_1(vw)v_1+u_1v(v_1w)+u_1w(v_1v)+u_1(v_1wv)+u_1(v_1(vw)) \
\end{align*}
and
\begin{align*}
-B&=u_1v_1wv+(u_1v_1)(vw) \\
&\equiv u_1wv_1v+u_1(v_1w)v+u_1v_1(vw) \ \\
&\equiv u_1wvv_1+u_1w(v_1v)+u_1v(v_1w)+u_1(v_1wv)+u_1v_1(vw). \ \
\
\end{align*}
So,
\begin{equation*}
(f_{uvw},f_{u_1v_1w_1})_{uvw}=A+B\equiv
u_1(vw)v_1+u_1(v_1(vw))-u_1v_1(vw)\equiv0 \ \ mod(S_0, uvw).
\end{equation*}

This completes our proof. \ \ \ \ $\square$

\end{document}